\documentclass[a4paper,12pt]{amsart}

\usepackage{amssymb,amscd,amsthm,mathrsfs}
\usepackage{epsfig,graphicx,psfrag,mathrsfs}
\usepackage[dvips]{color}
\usepackage[ansinew]{inputenc}
\usepackage[centertags]{amsmath}

 \def\RR{{\mathbb R}} \def\SS{{\mathbb S}}
\def\TT{{\mathbb T}}    
 \def\ZZ{{\mathbb Z}}

\def\cB{\mathcal{B}}

   \def\cR{\mathcal{R}}

\newtheorem*{teo*}{Theorem}
\newtheorem*{teoA}{Theorem A}
\newtheorem*{teoB}{Proposition B}
\newtheorem{teo}{Theorem}[section]
\newtheorem{ad}[teo]{Addendum}

\newtheorem{cor}[teo]{Corollary}
\newtheorem*{af}{Claim}
\newtheorem{lema}[teo]{Lemma}
\newtheorem{prop}[teo]{Proposition}

\newcommand{\bi}{\begin{itemize}}
\newcommand{\ei}{\end{itemize}}

\theoremstyle{definition}

\theoremstyle{remark}
\newtheorem{obs}[]{Remark}

\newcommand{\demo}[1]{\vspace{.05in}{\sc\noindent Proof #1.}}
\newcommand{\dem}{\vspace{.05in}{\sc\noindent Proof.}}
\newcommand{\lqqd}{\par\hfill {$\Box$} \vspace*{.05in}}
\newcommand{\finobs}{\par\hfill{$\diamondsuit$} \vspace*{.05in}}

\newcommand{\eps}{\varepsilon}

\newcommand{\en}{\subset}

\DeclareMathOperator{\Homeo}{Homeo}
\DeclareMathOperator{\modulo}{mod}

\author[R. Potrie]{Rafael Potrie}
\address{CMAT, Facultad de Ciencias, Universidad de la Rep\'ublica, Uruguay}\curraddr{LAGA; Institute Galilee, Universite Paris 13, Villetaneuse, France}
\email{rpotrie@cmat.edu.uy}

\title{Recurrence of non-resonant homeomorphisms on the torus}
\thanks{The autor was partially supported by ANR Blanc DynNonHyp BLAN08-2$\_$313375 and  ANII Proyecto FCE2007$\_$577}

\begin{document}

\maketitle

\begin{abstract}
We prove that a homeomorphism of the torus homotopic to the identity whose rotation set is reduced to a single totally irrational vector is chain-recurrent. In fact, we show that pseudo-orbits can be chosen with a small number of jumps, in particular, that the nonwandering set is weakly transitive. We give an example showing that the nonwandering set of such a homeomorphism may not be transitive.
\end{abstract}

\section{Introduction}
We consider $\Homeo_0(\TT^2)$ to be the set of homeomorphisms homotopic to the identity. We shall say that $f\in \Homeo_0(\TT^2)$ is \emph{non-resonant} if the rotation set of $f$ is a unique vector $(\alpha,\beta)$ and the values $1, \alpha, \beta$ are irrationally independent (i.e. $\alpha, \beta$ and $\alpha/\beta$ are not rational). This ammounts to say that given any lift $F$ of $f$ to $\RR^2$, for every $z\in \RR^2$ we have that:

\begin{equation}\label{equationRotacion}  \lim_{n\to \infty} \frac{F^n(z)-z}{n} = (\alpha, \beta) (\modulo \ZZ^2) \end{equation}

In general, one can define the rotation set of a homeomorphism homotopic to the identity (see \cite{MZ}). In fact, although we shall not make it explicit, our constructions work in the same way for homeomorphisms of the torus whose rotation set is contained in a segment of slope $(\alpha,\beta)$ with $\alpha, \beta$ and $\alpha/\beta$ irrational and not containing zero.

Non-resonant torus homeomorphisms\footnote{These are called \emph{irrational pseudo-rotations} by several authors, but since some of them use the term exclusively for conservative ones, we adopt the definition used in \cite{Kwakkel}.} have been intensively studied in the last years looking for resemblance between them and homeomorphisms of the circle with irrational rotation number (see \cite{Kwap}, \cite{LeCalvez}, \cite{Jager}) and also constructing examples showing some difference between them (see \cite{Fayad}, \cite{BCL}, \cite{BCJL}, \cite{Jager2}).

In \cite{Kwakkel} the possible topologies of minimal sets these homeomorphisms admit are classified and it is shown that under some conditions, these minimal sets are unique and coincide with the non-wandering set\footnote{A point $x$ is \emph{wandering} for a homeomorphism $f$ if there exists a neighborhood $U$ of $x$ such that $f^n(U) \cap U = \emptyset$ for every $n\neq 0$. The \emph{non-wandering} set is the closed set of points which are not wandering.}. However, there is one kind of topology of minimal sets where the question of the uniqueness of minimal sets remains unknown. When the topology of a minimal set is of this last kind, \cite{BCJL} constructed an example where the non wandering set does not coincide with the unique minimal set, in fact, they construct a transitive non-resonant torus homeomorphism containing a proper minimal set as a skew product over an irrational rotation.

A natural example of non-resonant torus homeomorphism is the one given by a homeomorphism semiconjugated to an irrational rotation by a continuous map homotopic to the identity. In \cite{Jager} it is proved that a non-resonant torus homeomorphism is semiconjugated to an irrational rotation under some quite mild hypothesis.

Under the hypothesis of being semiconjugated by a monotone map\footnote{A monotone map is a map whose preimages are all compact and connected.} which has points whose preimage is a singleton, it is not hard to show the uniqueness of a minimal set (see for example \cite{Kwakkel} Lemma 14). However, as shown by Roberts in \cite{Roberts}, a continuous monotone map may be very degenerate and thus even if there exist such a semiconjugation, it is not clear whether there should exist a unique minimal set nor the kind of recurrence the homeomorphisms should have. Moreover, for general non-resonant torus homeomorphisms, there does not exist a semiconjugacy to the irrational rotation (even when there is ``bounded mean motion'', see \cite{Jager2}).

Here, we give a simple and self-contained proof (based on some ideas of \cite{Kwakkel} but not on the classification of the topologies of the minimal sets) of a result which shows that even if there may be more than one minimal set, the dynamics is in some sense irreducible. Clearly, transitivity of $f$ may not hold for a general non-resonant torus homeomorphism (it may even have wandering points, as in the product of two Denjoy counterexamples; some more elaborate examples may be found in \cite{Kwakkel}), but we shall show that, in fact, these homeomorphisms are weakly transitive. For a homeomorphism $f$ we shall denote $\Omega(f)$ to the non-wandering set of $f$ (i.e. the set of points $x$ such that for every neighborhood $U$ of $x$ there exists $n>0$ with $f^n(U)\cap U \neq \emptyset$).

\begin{teoA}\label{MainTeo} Let $f\in \Homeo_0(\TT^2)$ be a non-resonant torus homeomorphism, then, $f|_{\Omega(f)}$ is weakly transitive.
\end{teoA}

Recall that for $h:M\to M$ a homeomorphism, and $K$ an $h-$invariant compact set, we say that $h|_K$ is \emph{weakly transitive} if given two open sets $U$ and $V$ of $M$ intersecting $K$, there exists $n>0$ such  that $h^n(U)\cap V \neq \emptyset$ (the difference with being transitive is that for transitivity one requires the open sets to be considered relative to $K$).

This allows to re-obtain Corollary E of \cite{Jager}:

\begin{cor}\label{CorolarioTransitividad} Let $f\in \Homeo_0(\TT^2)$ be a non-resonant torus homeomorphism such that $\Omega(f)=\TT^2$. Then, $f$ is transitive.
\end{cor}

In fact, as a consequence of weak-transitivity, we can obtain also the more well known concept of chain-transitivity for non-resonant torus homeomorphisms.

\begin{cor}\label{Addendum} Let $f\in \Homeo_0(\TT^2)$ be a non-resonant torus homeomorphism, then, $f$ is chain-transitive.
\end{cor}

Recall that a homeomorphism $h$ of a compact metric space $M$ is \emph{chain-transitive} if for every pair of points $x,y\in M$ and every $\eps>0$ there exists an $\eps-$pseudo-orbit $x= z_0, \ldots, z_n=y$ with $n\geq 1$ (i.e. $d(z_{i+1}, h(z_i))<\eps$).

\dem Consider two points $x,y\in M$ and $\eps>0$.

 We first assume that $x \neq y$ are both nonwandering points which shows the idea in a simpler way. From Theorem A we know that there exists a point $z$ and $n>0$ such that $d(z,f(x))<\eps$ and $d(f^{n+1}(z), y)<\eps$. We can then consider the $\eps-$pseudo-orbit: $\{x, z, \ldots, f^n(z), y\}$.

Now, for general $x, y \in \TT^2$ we consider $n_0\geq 1$ such that $d(f^{n_0+1}(x), \Omega(f))<\eps/2$ and $d(f^{-n_0}(y), \Omega(f))<\eps/2$. Now, by Theorem A there exists $z\in \TT^2$ and $n>0$ such that $d(z,f^{n_0}(x))<\eps$ and $d(f^{n+1}(z),f^{-n_0}(y))<\eps$. Considering the following $\eps-$psudo-orbit $\{x, \ldots, f^{n_0-1}(x), z, \ldots, f^n(z), f^{-n_0}(y), \ldots, y\}$ we obtain a pseudo-orbit from $x$ to $y$ and thus proving chain-transitivity.
\lqqd

\begin{obs} We have proved that in fact, for every $\eps>0$ the pseudo-orbit can be made with only two ``jumps''.
\end{obs}

As a consequence of our study, we obtain the following result which may be of independent interest:

\begin{teoB}\label{PropConexosCpctsSeCortan} Let $f \in \Homeo_0(\TT^2)$ be a non-resonant torus homeomorphism and $\Lambda_1$ a compact connected set  such that $f(\Lambda_1) \en \Lambda_1$. Then, for every $U$ connected neighborhood of $\Lambda_1$, there exists $K>0$ such that:
\bi
\item[-] If $\Lambda_2$ is a compact set which has a connected component in the universal cover of diameter larger than $K$ then\footnote{This holds  if $\Lambda_2$ is a connected set such that $f^i(\Lambda_2)\en \Lambda_2$ for some $i\in\ZZ$ for example.},
\ei \noindent $U \cap \Lambda_2 \neq \emptyset$.
\end{teoB}

One could wonder if the stronger property of $\Omega(f)$ being transitive may hold. However, in section \ref{SectionExample} we present an example  where $\Omega(f)$ is a Cantor set times $\SS^1$, but for which the nonwandering set is not transitive.

\textit{Acknowledgements:} I would like to thank Sylvain Crovisier and Martin Sambarino for their support and their important corrections to this text. I would also like to thank Tobias Jager who kindly exchanged communications related to this result. The referee carefully looked at the paper and gave suggestions to improve the presentation and clarity.

\section{Reduction of the proofs of Theorem A and Proposition B}

In this section we shall reduce the proofs of Theorem A and Proposition B to Proposition \ref{MainProp} and its Addendum \ref{AdeStronger}.

We shall use the word \emph{domain} to refer to an open and connected set. We shall say a domain $U \in \TT^2$ is \emph{inessential, simply essential or doubly essential} depending on whether the inclusion of $\pi_1(U)$ in $\pi_1(\TT^2)$ is isomorphic to $0, \ZZ$ or $\ZZ^2$ respectively\footnote{In \cite{Kwakkel} these concepts are called \emph{trivial, essential} and \emph{doubly-essential}.}. If $U$ is simply essential or doubly essential, we shall say it is \emph{essential}.

\begin{obs}\label{ObsDoublyEssentialsecortan} Notice that if $U$ and $V$ are two doubly essential domains, then $U\cap V\neq \emptyset$. This is because the intersection number of two closed curves is a homotopy invariant and given two non-homotopic curves in $\TT^2$, they have non-zero intersection number, thus, they must intersect. Since clearly, being doubly essential, $U$ and $V$ contain non homotopic curves, we get the desired result. \finobs
\end{obs}

We claim that Theorem \ref{MainTeo} can be reduced to the following proposition.

\begin{prop}\label{MainProp} Given $f\in \Homeo_0(\TT^2)$ a non-resonant torus homeomorphism and $U$ an open set such that $f(U) \en U$ and $U$ intersects $\Omega(f)$, then we have that $U$ has a connected component which is doubly essential.
\end{prop}

Almost the same proof also yields the following statement which will imply Proposition B:

\begin{ad}\label{AdeStronger} For $f$ as in Proposition \ref{MainProp}, if $\Lambda$ is a compact connected set such that $f(\Lambda) \en \Lambda$, then, for every connected open neighborhood $U$ of $\Lambda$, we have that $U$ is doubly-essential.
\end{ad}

Notice that the fact that $f(\Lambda)\en \Lambda$ for $\Lambda$ compact implies that it contains recurrent points, and in particular, $\Lambda\cap \Omega(f)\neq \emptyset$.

\demo{of Theorem A and Proposition B} Let us consider two open sets $U_1$ and $V_1$ intersecting $\Omega(f)$, and let $U = \bigcup_{n > 0} f^n(U_1)$ and $V = \bigcup_{n < 0} f^n(V_1)$. These sets verify that $f(U) \en U$ and $f^{-1}(V)\en V$ and both intersect the nonwandering set.

Proposition \ref{MainProp} (applied to $f$ and $f^{-1}$) implies that both $U$ and $V$ are doubly essential, so, they must intersect. This implies that for some $n > 0$ and $m < 0$ we have that $f^n(U_1) \cap f^m(V_1) \neq \emptyset$, so, we have that $f^{n-m}(U_1)\cap V_1 \neq \emptyset$ and thus $\Omega(f)$ is weakly transitive.

Proposition B follows directly from Addendum \ref{AdeStronger} since given a doubly-essential domain $U$ in $\TT^2$, there exists $K>0$ such that its lift $p^{-1}(U)$ intersects every connected set of diameter larger than $K$.
\lqqd

\begin{obs}\label{ObsDimensionmayor} Notice that in higher dimensions, Remark \ref{ObsDoublyEssentialsecortan} does not hold. In fact, it is easy to construct two open connected sets containing closed curves in every homotopy class which do not intersect. So, even if we could show a result similar to Proposition \ref{MainProp}, it would not imply the same result. \finobs
\end{obs}

\section{Proof of Proposition \ref{MainProp}}

Consider a non-resonant torus homeomorphism $f \in \Homeo_0(\TT^2)$, and let us assume that $U$ is an open set which verifies $f(U)\en U$ and $U \cap \Omega(f)\neq \emptyset$. % $d_0(\cdot, \cdot)$ the usual distance in $\RR^2$ and $d(\cdot, \cdot)$ the induced distance in $\TT^2$. %So, we get that $d(\partial U, f(\overline{U}))>\delta>0$.

Since $U\cap \Omega(f)\neq \emptyset$, for some $N>0$ we have that there is a connected component of $U$ which is $f^N$-invariant. We may thus assume from the start that $U$ is a domain such that $f(U)\en U$ and $U\cap \Omega(f) \neq \emptyset$.

Let $p:\RR^2 \to \TT^2$ be the canonical projection. Consider $U_0 \en p^{-1}(U)$ a connected component. We can choose $F$ a lift of $f$ such that $F(U_0) \en U_0$.

We shall denote $T_{p,q}$ to the translation by vector $(p,q)$, that is, the map from the plane such that $T_{p,q}(x)= x+(p,q)$ for every $x\in \RR^2$.

\begin{lema}\label{LemaNoTrivial} The domain $U$ is essential.
\end{lema}

\dem \ Consider $x\in U_0$ such that $p(x)\in \Omega(f)$. And consider a neighborhood $V\en U_0$ of $x$. Assume that there exists $n_0>0$ and $(p,q)\in \ZZ^2 \setminus \{(0,0)\}$ such that $F^{n_0}(V)\cap (V +(p,q)) \neq \emptyset$. Since $U_0$ is $F$-invariant, we obtain two points in $U_0$ which differ by an integer translation, and since $U_0$ is connected, this implies that $U$ contains a non-trivial curve in $\pi_1(\TT^2)$ and thus, it is essential.

To see that there exists such $n_0$ and $(p,q)$, notice that otherwise, since $x$ is not periodic (because $f$ is a non-resonant torus homeomorphism) we could consider a basis $V_n$ of neighborhoods of $p(x)$ such that $f^k(V_n) \cap V_n = \emptyset$ for every $0<  k \leq n$. Since $x$ is non-wandering, there exists some $k_n > n$ such that $f^{k_n}(V_n)\cap V_n \neq \emptyset$, but since we have that $F^{k_n}(V_n)\cap (V_n +(p,q)) = \emptyset$ for every $(p,q)\in \ZZ^2 \setminus \{(0,0)\}$, we should have that $F^{k_n}(V_n)\cap V_n \neq \emptyset$ for every $n$. Since $k_n \to \infty$, we get that $f$ has zero as rotation vector, a contradiction.
\lqqd

%\begin{obs} Clearly, this result holds equally for any domain containing a compact connected forward invariant set in its interior since we used only the fact that $U$ had points which where non-wandering and thus this Lemma works also in the hypothesis of Addendum \ref{AdeStronger} \finobs
%\end{obs}

We conclude the proof of by showing the following lemma which has some resemblance with Lemma 11 in \cite{Kwakkel}.

\begin{lema}\label{LemaNoEsencial} The domain $U$ is doubly-essential.
\end{lema}

\dem Assume by contradiction that $U$ is simply-essential.

Since the inclusion of $\pi_1(U)$ in $\pi_1(\TT^2)$ is non-trivial by the previous lemma, there exists a closed curve $\eta$ in $U$ such that when lifted to $\RR^2$ joins a point $x \in U_0$ with $x+(p,q)$ (which will also belong to $U_0$ because $\eta$ is contained in $U$ and $U_0$ is a connected component of $p^{-1}(U)$).

We claim that in fact, we can assume that $\eta$ is a simple closed curve and such that $g.c.d(p,q)=1$ (the greatest common divisor). In fact, since $U$ is open, we can assume that the curve we first considered is in general position, and by considering a subcurve, we get a simple one (maybe the point $x$ and the vector $(p,q)$ changed, but we shall consider the curve $\eta$ is the simple and closed curve from the start). Since it is simple, the fact that $g.c.d(p,q)=1$ is trivial.

If $\eta_0$ is the lift of $\eta$ which joins $x\in U_0$ with $x+(p,q)$, we have that it is compact, so, we get that

$$ \tilde \eta = \bigcup_{n\in \ZZ} T_{np,nq}\eta_0$$

\noindent is a proper embedding of $\RR$ in $\RR^2$. Notice that $\tilde \eta \en U_0$.

By extending to the one point compactification of $\RR^2$ we get by using Jordan's Theorem (see \cite{Moise} chapter 4) that $\tilde \eta$ separates $\RR^2$ in two disjoint unbounded connected components which we shall call $L$ and $R$ and such that their closures $L\cup \tilde \eta$ and $R\cup \tilde \eta$ are topologically a half plane (this holds by Sch\"onflies Theorem, see \cite{Moise} chapter 9).

Consider any pair $a,b$ such that\footnote{We accept division by $0$ as being infinity.} $\frac{a}{b} \neq \frac{p}{q}$, we claim that $T_{a, b} (\tilde \eta) \cap U_0 = \emptyset$. Otherwise, the union $T_{a,b}(\tilde \eta) \cup U_0$ would be a connected set contained in $p^{-1}(U)$ thus in $U_0$ and we could find a curve in $U_0$ joining $x$ to $x+(a,b)$ proving that $U$ is doubly essential (notice that the hypothesis on $(a,b)$ implies that $(a,b)$ and $(p,q)$ generate a subgroup isomorphic to $\ZZ^2$), a contradiction.

Translations are order preserving, this means that $T_{a,b}(R) \cap R$ and $T_{a,b}(L)\cap L$ are both non-empty and either $T_{a,b}(R)\en R$ or $T_{a,b}(L) \en L$ (both can only hold in the case $\frac{a}{b} = \frac{p}{q}$). Also, one can easily see that $T_{a,b}(R)\en R$ implies that $T_{-a,-b}(L)\en L$.

Now, we choose $(a,b)$ such that there exists a curve $\gamma$ from $x$ to $x+(a,b)$ satisfying:

\bi
\item[-] $T_{a,b}(\tilde \eta) \en L$.
\item[-] $\gamma$ is disjoint from $T_{p,q}(\gamma)$.
\item[-] $\gamma$ is disjoint from $T_{a,b}(\tilde \eta)$ and $\tilde \eta$ except at its boundary points.
\ei

We consider $\tilde \eta_1 = T_{a,b}(\tilde \eta)$ and $\tilde \eta_2 = T_{-a,-b}(\tilde \eta)$. Also, we shall denote $\tilde \gamma= \gamma \cup T_{-a,-b}(\gamma)$ which joins $x-(a,b)$ with $x+(a,b)$.

We obtain that $U_0$ is contained in $\Gamma= T_{a,b}(R)\cap T_{-a,-b}(L)$ a band whose boundary is $\tilde \eta_1 \cup \tilde \eta_2$.

Since $U_0$ is contained in $\Gamma$ and is $F$-invariant, for every point $x\in U_0$ we have that $F^n(x)$ is a sequence in $\Gamma$, and since $f$ is a non-resonant torus homeomorphism, we have that $\lim \frac{F^n(x)}{n} = \lim \frac{F^n(x)-x}{n}= (\alpha,\beta)$ is totally irrational.

However, we notice that $\Gamma$ can be written as:

$$ \Gamma = \bigcup_{n \in \ZZ} T_{np, nq} (\Gamma_0)$$

\noindent where $\Gamma_0$ is a compact set in $\RR^2$. Indeed, if we consider the curve $\tilde \gamma \cup T_{a,b}(\eta_0) \cup T_{p,q}(\tilde \gamma) \cup T_{-a,-b}(\eta_0)$ we have a Jordan curve. Considering $\Gamma_0$ as the closure of the bounded component we have the desired fundamental domain.

So, if we consider a sequence of points $x_n \in \Gamma$ such that $\lim \frac{x_n}{n}$ exists and is equal to $v$ it will verify that the coordinates of $v$ have the same proportion as $p/q$, thus cannot be totally irrational. This is a contradiction and concludes the proof of the Lemma.

\lqqd

%\begin{obs} Considering $\Lambda$ as in Addendum \ref{AdeStronger} we see that this proof works equally well since we only used that $U_0$ contained points which remained there to create the non wanted rotation vector and not that the whole of $U_0$ was invariant.
%\finobs
%\end{obs}

We conclude this section by showing how the proof adapts to the case stated in Addendum \ref{AdeStronger}. Consider a compact connected set $\Lambda$ such that $f(\Lambda) \en \Lambda$, then, we have that $\Lambda$ contains points which are recurrent\footnote{Since it is a compact invariant set, it contains a minimal set whose points will be all recurrent.}.

Let $\tilde \Lambda$ be a connected component of $p^{-1}(\Lambda)$ which is $F$-invariant. Now, if $U$ is an open connected neighborhood of $\Lambda$ and $U_0$ is a connected component of $p^{-1}(U)$ containing $\tilde \Lambda$. Notice that $d(\partial U, \Lambda) > \delta >0$ so $d(\partial U_0, \tilde \Lambda) > \delta$ also.

Now, the same argument in Lemma \ref{LemaNoTrivial} can be used in order to show that $U$ must be essential: We can choose a point $x\in \tilde \Lambda \en U_0$ such that $p(x)$ is recurrent and the same argument shows that there will exist $(p,q) \in \ZZ^2 \setminus \{(0,0)\}$ such that $F^{n_0}(x)$ is $\delta$-close to $x +(p,q)$ and since $F^{n_0}(x)$ must be contained in $\tilde \Lambda$ we get that $x + (p,q)$ is contained in $U_0$ showing that $U$ is essential.

The proof that in fact $U$ is doubly-essential is now the same as in Lemma \ref{LemaNoEsencial} since one can see that invariance of $U$ was not used in the proof, one only needs that there are points in $U_0$ such that the orbits by $F$ remain in $U_0$ and this holds for every point in $\tilde \Lambda$.

\section{An example where $f|_{\Omega(f)}$ is not transitive}\label{SectionExample}

The example is similar to the one in section 2 of \cite{Jager2}, however, we do not know a priori if our specific examples admit or not a semiconjugacy.

Consider $g_1:S^1 \to S^1$ and $g_2: S^1 \to S^1$ Denjoy counterexamples with rotation numbers $\rho_1$ and $\rho_2$ which are irrationally independent and have minimal invariant sets $M_1$ and $M_2$ properly contained in $S^1$. We shall consider the following skew-product map $f_\beta:\TT^2 \to \TT^2$ given by:

$$ f_\beta(s,t) = (g_1(s), \beta(s)(t)) $$

\noindent where $\beta: S^1 \to \Homeo_+(S^1)$ is continuous and such that $\beta(s)(t)=g_2(t)$ for every $(s,t)\in M_1 \times S^1$.

The same proof as in Lemma 2.1 of \cite{Jager2} yields:

\begin{lema}\label{LemaUnicoMinimal} The map $f_\beta$ is a non-resonant torus homeomorphism and $M_1 \times M_2$ is the unique minimal set.
\end{lema}

\dem The proof is the same as the one in Lemma 2.1 of \cite{Jager2}. Indeed any invariant measure for $f$ must be supported in $M_1 \times M_2$ and the dynamics there is the product of two Denjoy counterexamples and thus uniquely-ergodic. Since rotation vectors can be computed with ergodic measures, we also get that $f_\beta$ has a unique rotation vector $(\rho_1,\rho_2)$ which is totally irrational by hypothesis.
\lqqd

Clearly, if we restrict the dynamics of $f_\beta$ to $M_1 \times S^1$ it is not hard to see that the nonwandering set will be $M_1\times M_2$ (it is a product system there). So, we shall prove that if $\beta$ is properly chosen, we get that $\Omega(f_\beta)= M_1 \times S^1$. In fact, instead of constructing a specific example, we shall show that for ``generic'' $\beta$ in certain space, this is satisfied, this will give the existence of such a $\beta$.

First, we define $\cB$ to be the set of continuous maps $\beta: S^1 \to \Homeo_+(S^1)$ such that $\beta(s)=g_2$ for every $s\in M_1$. We endow $\cB$ with the topology given by restriction from the set of every continuous map from $S^1$ to $\Homeo_+(S^1)$. With this topology, $\cB$ is a closed subset of the set of continuous maps from $S^1 \to \Homeo_+(S^1)$ which is a Baire space, thus, $\cB$ is a Baire space.

So, the existence of the desired $\beta$ is a consequence of:

\begin{lema} There exists a dense $G_\delta$ (residual) subset of $\cB$ of maps such that the induced map $f_\beta$ verifies that $\Omega(f_\beta)= M_1 \times \SS^1$.
\end{lema}

\dem \ First, we will prove the lemma assuming the following claim:

\begin{af} Given $\beta \in \cB$, $x\in M_1 \times S^1$, $\eps>0$ and $\delta>0$ there exists $\beta' \in \cB$ which is $\delta-$close to $\beta$ such that there exists $k>0$ with $f_{\beta'}^k (B(x,\eps)) \cap B(x,\eps) \neq \emptyset$.
\end{af}

Assuming this claim, the proof of the Lemma is a standard Baire argument: Consider $\{x_n\} \en M_1\times S^1$ a countable dense set. Using the claim, we get that the sets $\cB_{n,N}$ consisting of the functions $\beta \in \cB$ such that there exists a point $y$ and a value $k>0$ such that $y$ and $f^k_{\beta}(y)$ belong to $B(x_n,1/N)$ is a dense set. Also, the set $\cB_{n,N}$ is open, since the property is clearly robust for $C^0$ perturbations of $f_\beta$. This implies that the set $\cR= \bigcap_{n,N} \cB_{n,N}$ is a residual set, which implies, by Baire's theorem that it is in fact dense.

For $\beta \in \cR$ we get that given a point $x\in M_1\times S^1$ and $\eps>0$, we can choose $x_n \in B(x,\eps/2)$ and $N$ such that $1/N<\eps/2$. Since $\beta \in \cB_{n,N}$ we have that there exists $k>0$ such that $f_\beta^k(B(x,\eps)) \cap B(x,\eps) \neq \emptyset$ proving that $M_1 \times S^1$ is nonwandering for $f_\beta$ as desired.

\demo{of the Claim} The point $x \in M_1 \times S^1$ can be written as $(s,t)$ in the canonical coordinates.

Choose an interval $(a,b) \en (s-\eps, s+\eps)$ contained in a wandering interval of $g_1$. Then, there exists a sequence of integers $k_n\to +\infty$ such that $g_1^{k_n}((a,b)) \en (s-\eps,s+\eps)$ for all $n\geq 0$. Further, the orbits of $a$ and $b$ are disjoint and do not belong to $M_1$. Let $\gamma= (a,b) \times \{t\}$.

We can assume that $f_\beta^{k_n}(\gamma) \cap B(x,\eps) = \emptyset$ for every $n>0$, otherwise, there is nothing to prove.

We shall thus consider a $\delta-$perturbation of $\beta$ such that it does not modify the orbit of $(a,t)$ but moves the orbit of $(b,t)$ in one direction making it give a complete turn around $\SS^1$ and thus an iterate of $\gamma$ will intersect $B(x,\eps)$.

Let $s_n = g_1^n(b)$ and $\beta^n(s_0) = \beta(s_{n-1})\circ \ldots \circ \beta(s_0)$. Note that $\frac{1}{n}\beta^n(s_0)(t) \to \rho(g_2)$ as $n\to \infty$ since $\beta(s_k) \to g_2$ as $k\to \infty$. At the same time, if we let

$$\beta_\theta^n(s_0)= R_\theta \circ \beta(s_{n-1}) \circ R_\theta \circ \beta_{s_{n-2}} \circ \ldots \circ R_\theta \circ \beta(s_0) $$

Then, $\beta_\theta^n(s_0)(t) \to \rho' > \rho(g_2)$ since $R_\theta \circ \beta_{s_k}$ converges to $R_\theta \circ g_2$ which has rotation number strictly greater than $g_2$ (see for example \cite{KH} Proposition 11.1.9). If we denote by $\tilde \beta^n$, respectively $\tilde \beta^n_\theta$ the lifts of $\beta^n$ and $\beta_\theta^n$ to $\RR$, then, this implies that there exists $n_0$ such that for $n>n_0$ one has

$$|\tilde \beta_\theta^n(t) - \tilde \beta^n(t) | > 1$$

So, if we consider $k_n>n_0$ and we choose $\beta'$ such that:

\bi
\item[-] it coincides with $\beta$ in the $g_1$-orbit of $a$,
\item[-] it coincides with $R_\theta \circ \beta$ in the points $\{b, g_1(b), \ldots, g_1^{k_n}(b)\}$,
\item[-] is at distance smaller than $\delta$ from $\beta$,
\ei

we have that $f^{k_n}_{\beta'}(\gamma) \cap B(x, \eps) \neq \emptyset$ as desired.

\finobs
\lqqd


\begin{thebibliography}{2}

\bibitem[BCJL]{BCJL} F.Beguin, S.Crovisier, T.Jager and F.Le Roux, Denjoy constructions for fibred homeomorphisms of the torus, \emph{Trans. Amer. Math. Soc.} {\bf 361} 11 (2009) 5851-5883.

\bibitem[BCL]{BCL} F. Beguin, S.Crovisier and F. Le Roux, Construction of curious minimal uniquely ergodic homeomorphisms on manifolds: the Denjoy-Rees technique. \emph{Ann. Sci. ENS.} {\bf 40} (2007), 251-308.

\bibitem[F]{Fayad} B. Fayad. Weak mixing for reparameterized linear flows on the torus. \emph{Ergodic Theory
Dyn. Syst.}, {\bf 22} 1 (2002), 187–201.

\bibitem[J1]{Jager} T. Jager, Linearisation of conservative toral homeomorphisms, \emph{Inventiones Math.} {\bf 176}(3) (2009), 601-616.

\bibitem[J2]{Jager2} T. Jager, The concept of bounded mean motion for toral homeomorphisms, \emph{Dynamical Systems. An International Journal} {\bf 24}3   (2009), 277-297.

\bibitem[KH]{KH} A. Katok and B. Hasselblatt, \emph{Introduction to the modern theory of dynamical systems}, Cambridge University Press (1995).

\bibitem[K1]{Kwakkel} F. Kwakkel, Minimal sets of non-resonant torus homeomorphisms, \emph{Fund. Math.} {\bf 211} (2011), 41-76.

\bibitem[K2]{Kwap} J. Kwapisz. Combinatorics of torus diffeomorphisms. \emph{Ergodic Theory Dynam. Systems} {\bf 23}
(2003), 559-586.

\bibitem[L]{LeCalvez} P. Le Calvez, . Ensembles invariants non enlaces des diffeomorphismes du tore et de l'anneau.
\emph{Invent. Math.} {\bf 155} (2004), 561-603.


\bibitem[MZ]{MZ} M. Misiurewicz and K. Ziemian. Rotation sets for maps of tori. \emph{J. Lond. Math. Soc.}, {\bf 40}(1989) 490–506.

\bibitem[M]{Moise} E. Moise, \emph{Geometric topology in dimensions 2 and 3}, Graduate texts in Mathematics {\bf 47} Springer (1977).


%\bibitem[M2]{Moore} R.L Moore, Concerning upper semi-continuous collections of continua, \emph{Trans. Amer. Math. Soc.} {\bf 27} (1925) 416-428.


\bibitem[R]{Roberts} J.H. Roberts, Collections filling a plane, \emph{Duke Math J.}  {\bf 2} n1 (1939) pp. 10-19.


\end{thebibliography}
\end{document}